\documentclass[leqno,english]{smfart}

\usepackage[leqno]{amsmath}
\usepackage{amssymb, amsthm, enumerate, latexsym, mathrsfs, mdwlist, stmaryrd}
\usepackage{ae,aecompl}
\usepackage[all]{xy}

\theoremstyle{remark}
\newtheorem{nul}{}[section]

\newtheorem*{rem*}{Remark}

\theoremstyle{definition}
\newtheorem{dfn}[nul]{Definition}

\newtheorem{ntn}[nul]{Notation}

\newtheorem*{dfn*}{Definition}
\newtheorem*{axm*}{Axiom}
\newtheorem*{ntn*}{Notation}
\newtheorem*{exm*}{Example}
\newtheorem*{exr*}{Exercise}
\newtheorem*{int*}{Intuition}
\newtheorem*{qst*}{Question}

\theoremstyle{plain}

\newtheorem{thm}[nul]{Theorem}
\newtheorem{prp}[nul]{Proposition}
\newtheorem{cor}[nul]{Corollary}
\newtheorem{lem}[nul]{Lemma}

\newtheorem*{thm*}{Theorem}
\newtheorem*{prp*}{Proposition}
\newtheorem*{cor*}{Corollary}
\newtheorem*{lem*}{Lemma}
\newtheorem*{cnj*}{Conjecture}

\numberwithin{equation}{nul}

\swapnumbers

\DeclareMathOperator{\cell}{cell}
\DeclareMathOperator{\cof}{cof}

\DeclareMathOperator{\fib}{fib}

\DeclareMathOperator{\Ho}{Ho}

\DeclareMathOperator{\holim}{holim}

\DeclareMathOperator{\inj}{inj}

\DeclareMathOperator{\Mor}{Mor}

\DeclareMathOperator{\Obj}{Obj}

\DeclareMathOperator{\RMor}{\mathbf{R}Mor}

\DeclareMathOperator{\Sect}{Sect}

\newcommand{\Prod}{\prod}

\newcommand{\coprd}{\amalg}

\newcommand{\CC}{\mathbf{C}}
\newcommand{\DD}{\mathbf{D}}
\newcommand{\EE}{\mathbf{E}}
\newcommand{\FF}{\mathbf{F}}

\newcommand{\LL}{\mathbf{L}}
\newcommand{\MM}{\mathbf{M}}
\newcommand{\NN}{\mathbf{N}}

\newcommand{\RR}{\mathbf{R}}

\newcommand{\XX}{\mathbf{X}}

\newcommand{\op}{\mathrm{op}}

\newcommand{\fromto}[2]{\xymatrix@1@C=18pt{{#1}\ar[r]&{#2}}}
\newcommand{\goesto}[2]{\xymatrix@1@C=12pt{{#1}\,\ar@{|->}[r]&{#2}}}
\newcommand{\adjunct}[4]{{#1}:\xymatrix@1@C=18pt{{#2}\ar@<0.5ex>[r]&{#3}\ar@<0.5ex>[l]}:{#4}}

\begin{document}

\title{On the Dreaded Right Bousfield Localization}
\author{Clark Barwick}
\address{Matematisk Institutt\\
Universitetet i Oslo\\
Boks 1053 Blindern\\
0316 Oslo\\
Norge}
\curraddr{School of Mathematics\\
Institute for Advanced Study\\
Einstein Drive\\
Princeton, NJ 08540-0631\\
USA}
\email{clarkbar@gmail.com}
\thanks{This work was supported by a research grant from the Yngre fremragende forskere, administered by J. Rognes at the Matematisk Institutt, Universitetet i Oslo.}

\begin{abstract} I verify the existence of right Bousfield localizations of right semimodel categories, and I apply this to construct a model of the homotopy limit of a left Quillen presheaf as a right semimodel category.
\end{abstract}

\maketitle
\thispagestyle{empty}

The right Bousfield localization --- or colocalization --- of a model category $\MM$ with respect to a set $K$ of objects is a model for the homotopy theory generated by $K$ --- i.e., of objects that can be written as a homotopy colimit of objects of $K$. Unfortunately, the right Bousfield localization need not exist as a model category unless $\MM$ is right proper. This is a rather severe limitation, as some operations on model categories --- such as left Bousfield localization --- tend to destroy right properness, and as many interesting model categories are not right proper.

In this very brief note, I show that the right Bousfield localization of a model category $\MM$ naturally exists instead as a right semimodel category --- i.e., as categories with classes of weak equivalences, cofibrations, and fibrations that satisfy all the usual axioms of a model category apart from the existence of factorizations into cofibrations followed by trivial fibrations and the lifting criterion for trivial fibrations, unless the targets of the trivial fibrations are fibrant. This result holds with no properness assumptions on $\MM$.

As an application of the right Bousfield localization, I produce a good model for the homotopy limit of left Quillen presheaves.

Thanks to J. Bergner, P. A. {\O}stv{\ae}r, and B. Toën for persistent encouragement and hours of interesting discussion. Thanks especially to M. Spitzweck for a profound and lasting impact on my work; were it not for his insights and questions, there would be nothing for me to report here or anywhere else.

\setcounter{tocdepth}{2}
\tableofcontents

\section{A taxonomy of homotopy theory} It is necessary to establish some general terminology for categories with weak equivalences and various bits of extra structure. This terminology includes such arcane and baroque concepts as structured homotopical categories and semimodel categories.

\setcounter{nul}{-1}

\begin{nul} Suppose $\XX$ a universe.
\end{nul}

\subsection*{Structured homotopical categories}{Here I define the general notion of structured homotopical categories. Structured homotopical categories contain lluf subcategories of cofibrations, fibrations, and weak equivalences, satisfying the easiest conditions on model categories.}

\begin{dfn}\label{dfn:strhomotopical} Suppose $(\CC,w\CC)$ a homotopical $\XX$-category equipped with two lluf subcategories $\cof\CC$ and $\fib\CC$.
\begin{enumerate}[(\ref{dfn:strhomotopical}.1)]
\item Morphisms of $\cof\CC$ (respectively, of $\fib\CC$) are called \emph{cofibrations} (resp., \emph{fibrations}).
\item Morphisms of $w\cof\CC:=w\CC\cap\cof\CC$ (respectively, of $w\fib\CC:=w\CC\cap\fib\CC)$ are called \emph{trivial cofibrations} (resp., \emph{trivial fibrations}).
\item Objects $X$ of $\CC$ such that the morphism $\fromto{\varnothing}{X}$ (respectively, the morphism $\fromto{X}{\star}$) is an element of $\cof\CC$ (resp., of $\fib\CC$) are called \emph{cofibrant} (resp. \emph{fibrant}); the full subcategory comprised of all such objects will be denoted $\CC_c$ (resp., $\CC_f$).
\item In the context of a functor $\fromto{\EE}{\CC}$, a morphism (respectively, an object) of $\EE$ will be called a \emph{$\CC$-weak equivalence}, a \emph{$\CC$-cofibration}, or a \emph{$\CC$-fibration} (resp., \emph{$\CC$-cofibrant} or \emph{$\CC$-fibrant}) if its image under $\fromto{\EE}{\CC}$ is a weak equivalence, a cofibration, or a fibration (resp., cofibrant or fibrant) in $\CC$, respectively. The full subcategory of $\EE$ comprised of all $\CC$-cofibrant (respectively, $\CC$-fibrant) objects will be denoted $\EE_{\CC,c}$ (resp., $\EE_{\CC,f}$).
\item\label{item:strhomotopical} One says that $(\CC,w\CC,\cof\CC,\fib\CC)$ is a \emph{structured homotopical $\XX$-category} if the following axioms hold.
\begin{enumerate}[(\ref{dfn:strhomotopical}.\ref{item:strhomotopical}.1)]
\item The category $\CC$ contains all limits and colimits.
\item The subcategories $\cof\CC$ and $\fib\CC$ are closed under retracts.
\item The set $\cof\CC$ is closed under pushouts by arbitrary morphisms; the set $\fib\CC$ is closed under pullbacks by arbitrary morphisms.
\end{enumerate}
\end{enumerate}
\end{dfn}

\begin{lem} The data $(\CC,w\CC,\cof\CC,\fib\CC)$ is a structured homotopical $\XX$-category if and only if the data $(\CC^{\op},w(\CC^{\op}),\cof(\CC^{\op}),\fib(\CC^{\op}))$ is as well, wherein
\begin{equation*}
w(\CC^{\op}):=(w\CC)^{\op}\quad\cof(\CC^{\op}):=(\fib\CC)^{\op}\quad\fib(\CC^{\op}):=(\cof\CC)^{\op}.
\end{equation*}
\end{lem}

\begin{nul} One commonly refers to $\CC$ alone as a structured homotopical category, omitting the explicit reference to the data of $w\CC$, $\cof\CC$, and $\fib\CC$.
\end{nul}

\begin{ntn}\label{ntn_subcCpCpsubf} Suppose $\CC$ a structured homotopical category, $p$ an object of $\Delta$.
\begin{enumerate}[(\ref{ntn_subcCpCpsubf}.1)]
\item Denote by ${}_c(\CC^p)$ (respectively, $(\CC^p)_f$) the full subcategory of $\CC^p$ comprised of those diagrams $X:\fromto{p}{\CC}$ such that $X(0)$ is cofibrant (resp., such that $X(p)$ is fibrant).
\item In the context of a functor $\fromto{\EE}{\CC}$, denote by ${}_{\CC,c}(\EE^p)$ (respectively, $(\EE^p)_{\CC,f}$ the full subcategory of $\EE^p$ comprised of those diagrams $X:\fromto{p}{\EE}$ such that $X(0)$ is $\CC$-cofibrant (resp., such that $X(p)$ is $\CC$-fibrant).
\end{enumerate}
\end{ntn}

\subsection*{Semimodel categories}{Semimodel categories are structured homotopical categories that, like model categories, include lifting and factorization axioms, but only for particular morphisms. I now turn to a sequence of standard results from the homotopy theory of model categories, suitably altered to apply to semimodel categories. I learned nearly all of the following ideas and results from M. Spitzweck and his thesis \cite{math.AT/0101102}.}

\begin{dfn}\label{dfn:semimodel} Suppose $\CC$ and $\EE$ two structured homotopical $\XX$-categories.
\begin{enumerate}[(\ref{dfn:semimodel}.1)]
\item\label{item:leftDsemimodel} An adjunction
\begin{equation*}
F_{\CC}:\xymatrix@1@C=18pt{\EE\ar@<0.5ex>[r]&\CC\ar@<0.5ex>[l]}:U_{\CC}
\end{equation*}
is a \emph{left $\EE$-semimodel $\XX$-category} if the following axioms hold.
\begin{enumerate}[(\ref{dfn:semimodel}.\ref{item:leftDsemimodel}.1)]
\item The right adjoint $U_{\CC}$ preserves fibrations and trivial fibrations.
\item Any cofibration of $\CC$ with $\EE$-cofibrant domain is an $\EE$-cofibration.
\item The initial object $\varnothing$ of $\CC$ is $\EE$-cofibrant.
\item In $\CC$, any cofibration has the left lifting property with respect to any trivial fibration, and any fibration has the right lifting property with respect to any trivial cofibration with $\EE$-cofibrant domain.
\item There exist functorial factorizations of any morphism of $\CC$ into a cofibration followed by a trivial fibration and of any morphism of $\CC$ with $\EE$-cofibrant domain into a trivial cofibration followed by a fibration.
\end{enumerate}
\item\label{item:rightDsemimodel} An adjunction
\begin{equation*}
F_{\CC}:\xymatrix@1@C=18pt{\CC\ar@<0.5ex>[r]&\EE\ar@<0.5ex>[l]}:U_{\CC}
\end{equation*}
is a \emph{right $\EE$-semimodel $\XX$-category} if the corresponding adjunction
\begin{equation*}
U_{\CC}^{\op}:\xymatrix@1@C=18pt{\EE^{\op}\ar@<0.5ex>[r]&\CC^{\op}\ar@<0.5ex>[l]}:F_{\CC}^{\op}
\end{equation*}
is a left $\EE$-semimodel $\XX$-category.
\item\label{item:leftsemimodel} One says that $\CC$ is a(n) \emph{(absolute) left semimodel $\XX$-category} if the identity adjunction
\begin{equation*}
\xymatrix@1@C=18pt{\CC\ar@<0.5ex>[r]&\CC\ar@<0.5ex>[l]}
\end{equation*}
is a left $\CC$-semimodel category.
\item\label{item:rightsemimodel} One says that $\CC$ is a(n) \emph{(absolute) right semimodel $\XX$-category} if $\CC^{\op}$ is a left semimodel $\XX$-category.
\item One says that $\CC$ is a \emph{model $\XX$-category} if it is both a left and right semimodel $\XX$-category.
\end{enumerate}
\end{dfn}

\begin{nul} By the usual abuse, I refer to $\CC$ alone as the $\EE$-left semimodel or $\EE$-right semimodel $\XX$-category, omitting explicit mention of the adjunction.
\end{nul}

\begin{lem} Consider a pair of composable adjunctions of structure homotopical $\XX$-categories
\begin{equation*}
\xymatrix@1@C=18pt{\EE\ar@<0.5ex>[r]&\DD\ar@<0.5ex>[r]\ar@<0.5ex>[l]&\CC\ar@<0.5ex>[l]}\textrm{\quad (respectively,\quad}\xymatrix@1@C=18pt{\CC\ar@<0.5ex>[r]&\DD\ar@<0.5ex>[r]\ar@<0.5ex>[l]&\EE\ar@<0.5ex>[l]}\textrm{\quad);}
\end{equation*}
then if under these adjunctions $\CC$ is a left (resp., right) $\EE$-semimodel category, and $\DD$ is a left (resp., right) $\EE$-semimodel category, then $\CC$ is a left (resp., right) $\DD$-semimodel category under the composite.
\begin{proof} By duality it suffices to prove the statement for left semimodel categories, and since any $\DD$-cofibrant object of $\CC$ is also $\EE$-cofibrant, the result follows immediately.
\end{proof}
\end{lem}

\begin{lem}\label{lem:modeliffstarsemimod} The following are equivalent for a structured homotopical category $\CC$.
\begin{enumerate}[(\ref{lem:modeliffstarsemimod}.1)]
\item $\CC$ is a left $\star$-semimodel category.
\item $\CC$ is a right $\star$-semimodel category.
\item $\CC$ is a model category.
\end{enumerate}
\begin{proof} To be a left $\star$-semimodel category is exactly to have the lifting and factorization axioms with no conditions on the source of the morphism, hence to be a right semimodel category as well. The dual assertion follows as usual.
\end{proof}
\end{lem}

\begin{lem}\label{lem:liftingcharsemi} Suppose $\EE$ a structured homotopical category, $\CC$ a left (respectively, right) $\EE$-semimodel $\XX$-category.
\begin{enumerate}[(\ref{lem:liftingcharsemi}.1)]
\item A morphism $i:\fromto{K}{L}$ (resp., a morphism $f:\fromto{K}{L}$ with $\EE$-fibrant codomain $L$) has the left lifting property with respect to every trivial fibration (resp., every trivial fibration with $\EE$-fibrant codomain) if and only if $i$ is a cofibration.
\item Any morphism $i:\fromto{K}{L}$ with $\EE$-cofibrant domain $K$ (resp., any morphism $f:\fromto{K}{L}$) has the left lifting property with respect to every fibration if and only if $i$ is a trivial cofibration.
\item Any morphism $p:\fromto{Y}{X}$ with $\EE$-cofibrant domain $Y$ (resp., any morphism $p:\fromto{Y}{X}$) has the right lifting property with respect to every trivial cofibration with $\EE$-cofibrant domain (resp., every trivial cofibration) if and only if $p$ is a fibration.
\item A morphism $p:\fromto{Y}{X}$ (resp., a morphism $p:\fromto{Y}{X}$ with $\EE$-fibrant codomain) has the right lifting property with respect to every cofibration if and only if $p$ is a trivial fibration.
\end{enumerate}
\begin{proof} This follows immediately from the appropriate factorization axioms along with the retract argument.
\end{proof}
\end{lem}

\begin{cor}\label{cor:liftingcharsemi}
\begin{enumerate}[(\ref{cor:liftingcharsemi}.1)]
\item If $\CC$ is a left semimodel $\XX$-category, a morphism $p:\fromto{Y}{X}$ satisfies the right lifting property with respect to the trivial cofibrations with cofibrant domains if and only if there exists a trivial fibration $\fromto{Y'}{Y}$ such that the composite morphism $\fromto{Y'}{X}$ is a fibration.
\item Dually, if $\CC$ is a right semimodel $\XX$-category, a morphism $i:\fromto{K}{L}$ satisfies the left lifting property with respect to the trivial fibrations with fibrant codomains if and only if there exists a trivial cofibration $\fromto{L}{L'}$ such that the composite morphism $\fromto{K}{L'}$ is a cofibration.
\end{enumerate}
\begin{proof} The assertions are dual, so it is enough to prove the first. Morphisms satisfying a right lifting property are of course closed under composition. Conversely, suppose $\fromto{Y}{X}$ a morphism, $\fromto{Y'}{Y}$ a trivial fibration such that the composition $\fromto{Y'}{X}$ satisfies the left lifting property with respect to a trivial cofibration $\fromto{K}{L}$ is a trivial cofibration with cofibrant domain $K$. Then for any diagram
\begin{equation*}
\xymatrix@C=18pt@R=18pt{
K\ar[d]\ar[r]&Y\ar[d]\\
L\ar[r]&X,
}
\end{equation*}
there is a lift to a diagram
\begin{equation*}
\xymatrix@C=18pt@R=18pt{
&Y'\ar[d]\\
K\ar[ur]\ar[d]\ar[r]&Y\ar[d]\\
L\ar[r]&X.
}
\end{equation*}
By assumption there is a lift of the exterior quadrilateral, and this provides a lift of the interior square as well.
\end{proof}
\end{cor}

\begin{prp}[\protect{\cite[Proposition 2.4]{math.AT/0101102}}]\label{prp_homotopycatsemimodel} Suppose $\EE$ a left semimodel $\XX$-category, and suppose $\CC$ a left $\EE$-semimodel $\XX$-category. Suppose $f,g:\fromto{B}{X}$ two maps in $\CC$.
\begin{enumerate}[(\ref{prp_homotopycatsemimodel}.1)]
\item Suppose $f\overset{l}{\sim}g$; then $h\circ f\overset{l}{\sim}h\circ g$ for any morphism $h:\fromto{X}{Y}$ of $\CC$.
\item Dually, suppose $f\overset{r}{\sim}g$; then $f\circ k\overset{r}{\sim}g\circ k$ for any morphism $k:\fromto{A}{B}$ of $\CC$.
\item Suppose $B$ cofibrant, and suppose $h:\fromto{X}{Y}$ any morphism of $\CC_{\EE,c}$; then $f\overset{r}{\sim}g$ only if $h\circ f\overset{r}{\sim}h\circ g$.
\item Dually, suppose $X$ fibrant, and suppose $k:\fromto{A}{B}$ any morphism of $\CC_{\EE,c}$; then $f\overset{l}{\sim}g$ only if $f\circ h\overset{l}{\sim}g\circ h$.
\item If $B$ is cofibrant, then left homotopy is an equivalence relation on $\Mor(B,X)$.
\item If $B$ is cofibrant and $X$ is $\EE$-cofibrant, then $f\overset{l}{\sim}g$ only if $f\overset{r}{\sim}g$.
\item Dually, if $X$ is fibrant and $B$ is cofibrant $\EE$-cofibrant, then $f\overset{r}{\sim}g$ only if $f\overset{l}{\sim}g$.
\item If $B$ is cofibrant, $X$ is cofibrant $\EE$-cofibrant, and $h:\fromto{X}{Y}$ is either a trivial fibration or a weak equivalence between fibrant objects, then $h$ induces a bijection
\begin{equation*}
\fromto{(\Mor(B,X)/\overset{l}{\sim})}{(\Mor(B,Y)/\overset{l}{\sim})}.
\end{equation*}
\item Dually, if $A$ is $\EE$-cofibrant, $X$ is fibrant and $\EE$-cofibrant, and $k:\fromto{A}{B}$ is either a trivial cofibration with $A$ $\DD$-cofibrant or a weak equivalence between cofibrant objects, then $k$ induces a bijection
\begin{equation*}
\fromto{(\Mor(B,X)/\overset{r}{\sim})}{(\Mor(A,X)/\overset{r}{\sim})}.
\end{equation*}
\end{enumerate}

The obvious dual statements for right semimodel and $\EE$-right semimodel categories also hold.
\end{prp}

\begin{cor} Suppose $\EE$ a left (respectively, right) semimodel $\XX$-category, $\CC$ a left $\EE$-semimodel (resp., right $\EE$-semimodel) $\XX$-category. Then $\Ho\CC$ is an $\XX$-category.
\end{cor}

\begin{dfn} Suppose $\EE$ a left (respectively, right) semimodel $\XX$-category, $\CC$ a left $\EE$-semimodel (resp., right $\EE$-semimodel) $\XX$-category; suppose $h$ a homotopy class of morphisms of $\CC$. Then a morphism $f$ of $\CC$ is said to be a \emph{representative} of $h$ if the images of $f$ and $h$ are isomorphic as objects of the arrow category $(\Ho\MM)^1$.
\end{dfn}

\begin{dfn}\label{dfn_Quillenadjsemi} Suppose $\EE$ a structured homotopical category.
\begin{enumerate}[(\ref{dfn_Quillenadjsemi}.1)]
\item Suppose $\CC$ and $\CC'$ two left (respectively, right) $\EE$-semimodel $\XX$-categories. Then an adjunction
\begin{equation*}
F:\xymatrix@1@C=18pt{\CC\ar@<0.5ex>[r]&\CC'\ar@<0.5ex>[l]}:U
\end{equation*}
is a \emph{Quillen adjunction} if $U$ preserves fibrations and trivial fibrations (resp., if $F$ preserves cofibrations and trivial cofibrations).
\item Suppose $\CC$ and $\CC'$ two left or right $\EE$-semimodel categories, and suppose
\begin{equation*}
F:\xymatrix@1@C=18pt{\CC\ar@<0.5ex>[r]&\CC'\ar@<0.5ex>[l]}:U
\end{equation*}
a Quillen adjunction. Then the \emph{left derived functor} $\LL F$ of $F$ is the right\footnote{sic, unfortunately.} Kan extension of the composite
\begin{equation*}
\xymatrix@C=18pt{
\CC\ar[r]^F&\CC'\ar[r]&\Ho\CC'
}
\end{equation*}
along the functor $\fromto{\CC}{\Ho\CC}$, and, dually, the \emph{right derived functor} $\RR U$ of $U$ is the left Kan extension of the composite
\begin{equation*}
\xymatrix@C=19pt{
\CC'\ar[r]^U&\CC\ar[r]&\Ho\CC
}
\end{equation*}
along the functor $\fromto{\CC'}{\Ho\CC'}$.
\end{enumerate}
\end{dfn}

\begin{prp}[\protect{\cite[p. 12]{math.AT/0101102}}] Suppose $\EE$ a structured homotopical category, and suppose $\CC$ and $\CC'$ two left or right semimodel $\EE$-categories, and suppose
\begin{equation*}
F:\xymatrix@1@C=18pt{\CC\ar@<0.5ex>[r]&\CC'\ar@<0.5ex>[l]}:U
\end{equation*}
a Quillen adjunction. Then the derived functors
\begin{equation*}
\LL F:\xymatrix@1@C=18pt{\Ho\CC\ar@<0.5ex>[r]&\Ho\CC'\ar@<0.5ex>[l]}:\RR U
\end{equation*}
exist and form an adjunction.
\end{prp}

\begin{dfn}\label{dfn_cofgensemi} Suppose $\EE$ a homotopical $\XX$-category, and suppose $\CC$ a left (respectively, right) $\EE$-semimodel $\XX$-category. Suppose, in addition, that $\lambda$ is a regular $\XX$-small cardinal.
\begin{enumerate}[(\ref{dfn_cofgensemi}.1)]
\item\label{item_cofgensemi}One says that $\CC$ is \emph{$\lambda$-tractable} if the underlying $\XX$-category of $\CC$ is locally $\lambda$-presentable, and if there exist $\XX$-small sets $I$ and $J$ of morphisms of $\CC_{\lambda}\cap\CC_{\EE,c}$ (resp., of $\CC_{\lambda}$) such that the following hold.
\begin{enumerate}[(\ref{dfn_cofgensemi}.\ref{item_cofgensemi}.1)]
\item A morphism (resp., a morphism with $\EE$-fibrant codomain) satisfies the right lifting property with respect to $I$ if and only if it is a trivial fibration.
\item A morphism satisfies the right lifting property with respect to $J$ if and only if it is a fibration.
\end{enumerate}
\item An $\XX$-small full subcategory $\CC_0$ of $\CC$ is \emph{homotopy $\lambda$-generating} if every object of $\CC$ is weakly equivalent to a $\lambda$-filtered homotopy colimit of objects of $\CC_0$.
\end{enumerate}

One says that $\CC$ is \emph{$\XX$-tractable} just in case there exists a regular $\XX$-small cardinal $\lambda$ for which it is $\lambda$-tractable. Likewise, an $\XX$-small full subcategory $\CC_0$ of $\CC$ is \emph{homotopy $\XX$-generating} if and only if for some sound $\XX$-doctrine $\lambda$, $\CC_0$ is homotopy $\lambda$-generating.
\end{dfn}

\begin{nul} If $\CC$ is a model category, then $\CC$ is $\lambda$-tractable \cite[1.3.1]{arXiv:0708.2067v1} if and only if the underlying \emph{left} $\EE$-semimodel category of $\CC$ is so. Observe, however, that a model category whose underlying right semimodel category is $\XX$-tractable need \emph{not} be $\XX$-tractable.
\end{nul}

\subsection*{The Reedy semimodel structure} If $A$ is a Reedy category and $\MM$ a left (respectively, right) semimodel category, the category $\MM^A$ of diagrams $\fromto{A}{\MM}$ has a Reedy left (resp., right) semimodel structure. The proofs of the results below are similar to the proofs of the classical results for model categories, save only that one must periodically insert the phrases ``with cofibrant domain'' (resp., ``with fibrant codomain''). I will therefore leave the proofs as an exercise.\footnote{Alternatively, see \cite[Propositions 2.6 and 2.7]{math.AT/0101102}, where the case of left semimodel categories is addressed.}

\begin{nul} Suppose $\XX$ a universe, $\MM$ a left (respectively, right) semimodel $\XX$-category.
\end{nul}

\begin{thm}\label{thm:projdirectinjinverse} 
\begin{enumerate}[(\ref{thm:projdirectinjinverse}.1)]
\item For any direct category $A$ \cite[1.1.2]{arxiv:0708.2832v1}, the functor category $\MM^A$ has its projective left (resp., right) semimodel structure, in which the weak equivalences and fibrations are defined objectwise.
\item For any inverse category $A$ \cite[1.1.3]{arxiv:0708.2832v1}, the functor category $\MM^A$ has its injective left (resp., right) semimodel structure, in which the weak equivalences and cofibrations are defined objectwise.
\end{enumerate}
\end{thm}

\begin{prp}\label{prp:leftrightKan} Suppose $f:\fromto{A}{B}$ a functor of $\XX$-small categories.
\begin{enumerate}[(\ref{prp:leftrightKan}.1)]
\item If $A$ and $B$ are direct categories, then the adjunction
\begin{equation*}
\adjunct{f_!}{\MM^A}{\MM^B}{f^{\star}}
\end{equation*}
is a Quillen adjunction between the projective left (resp., right) semimodel categories.
\item If $A$ and $B$ are inverse categories, then the adjunction
\begin{equation*}
\adjunct{f^{\star}}{\MM^B}{\MM^A}{f_{\star}}
\end{equation*}
is a Quillen adjunction between the injective left (resp., right) semimodel categories.
\end{enumerate}
\end{prp}

\begin{prp}\label{prp:cofsprojfibsinj}
\begin{enumerate}[(\ref{prp:cofsprojfibsinj}.1)]
\item For any direct category $A$, a morphism $\fromto{X}{Y}$ of the functor category $\MM^A$ is a cofibration or trivial cofibration in the projective left (resp., right) semimodel structure if and only if for any object $\alpha$ of $A$, the induced morphism
\begin{equation*}
\fromto{X_{\alpha}\coprd^{L_{\alpha}X}L_{\alpha}Y}{Y_{\alpha}}
\end{equation*}
is so.
\item For any inverse category $A$, a morphism $\fromto{X}{Y}$ of the functor category $\MM^A$ is a fibration or trivial fibration in the injective left (resp., right) semimodel structure if and only if for any object $\alpha$ of $A$, the induced morphism
\begin{equation*}
\fromto{X_{\alpha}}{M^{\alpha}X\times_{M^{\alpha}Y}Y_{\alpha}}
\end{equation*}
is so.
\end{enumerate}
\end{prp}

\begin{thm} Suppose $A$ an $\XX$-small Reedy category \cite[Definition 1.6]{arxiv:0708.2832v1}. Then the diagram category $\MM^A$ has its \emph{Reedy left (resp., right) semimodel structure}, in which a morphism $\phi:\fromto{X}{Y}$ is a weak equivalence, cofibration, or fibration if and only if both $i^{\to,\star}\phi$ in $M^{A^{\to}}$ and $i^{\gets,\star}\phi$ in $M^{A^{\gets}}$ are so.
\end{thm}

\section{The dreaded right Bousfield localization} The right Bousfield localization of a model category $\MM$ relative to a set of objects $K$ is ordinarily defined as a model category $R_K\MM$ equipped with a right Quillen functor $\fromto{\MM}{R_K\MM}$ satisfying a universal property dual to that of left Bousfield localizations.

\subsection*{Existence theorem} Suppose $\XX$ a universe, $\MM$ an $\XX$-cofibrantly generated model category, and $K$ an $\XX$-small set of objects of $\MM$. There does not appear to be an existence theorem for $R_K\MM$ unless $\MM$ is right proper. If $i:\fromto{A}{B}$ is a $K$-colocal cofibration of and $p:\fromto{Y}{X}$ is a $K$-colocal trivial fibration, it is necessary to show that for any diagram
\begin{equation*}
\xymatrix@C=18pt@R=18pt{
A\ar[d]\ar[r]&Y\ar[d]\\
B\ar[r]&X,
}
\end{equation*}
there exists a lift $\fromto{B}{Y}$. It turns out that this is easy to verify in case $X$ (and hence also $Y$) is fibrant in $\MM$. If $\MM$ is right proper, this is sufficient: $i$ has the left lifting property with respect to $p$ if and only if it has the left lifting property with respect to a replacement fibration $p':\fromto{Y'}{X'}$ of $p$ with $Y'$ and $X'$ fibrant \cite[Propositions 5.2.5 and 13.2.1]{MR2003j:18018}.

This leads one to the following observation: if one only seeks the left lifting property of $K$-colocal cofibrations with respect to $K$-colocal trivial fibrations \emph{with fibrant codomain}, then the right properness of $\MM$ is unnecessary here.

Likewise, the small object argument immediately provides factorizations into cofibrations followed by trivial fibrations when the codomain is fibrant. It is the existence of such factorizations for any morphism that requires right properness \cite[Proposition 5.3.5]{MR2003j:18018}.

Upon inspection of the standard proofs of the existence of $R_\MM$ for $\MM$ right proper, one can confirm that these are the only places where right properness is used. Hence $R_K\MM$ exists as a right $\MM$-semimodel category, even if $\MM$ is not right proper.

Thus, a simple modification of the traditional proof shows that $R_K\MM$ exists as a right $\MM$-semimodel category for any $\XX$-cofibrantly generated model category $\MM$ (or in fact for any model $\XX$-category $\MM$ satisfying Christensen and Isaksen's weaker condition \cite[Hypothesis 2.4]{MR2100680}), and that, as a right $\MM$-semimodel category, $R_K\MM$ is $\XX$-cofibrantly generated as well.\footnote{If $R_K\MM$ happens to be a model category, it does not seem to follow that $R_K\MM$ will be cofibrantly generated as a model category, unless some very strong conditions on $\MM$ are satisfied, e.g., that every object of $\MM$ be fibrant.} Here I give a complete proof of the existence of $R_K\MM$ as an $\XX$-tractable right $\CC$-semimodel category for any $\XX$-tractable right semimodel category $\CC$.

\begin{nul} Suppose $\EE$ a right semimodel $\XX$-category, $\CC$ an $\XX$-tractable right $\EE$-semimodel $\XX$-category, and $K$ a set of isomorphism classes of objects of $\Ho\CC$.
\end{nul}

\begin{dfn}\label{dfn:rBous}
\begin{enumerate}[(\ref{dfn:rBous}.1)]
\item If $H$ is a set of homotopy classes of morphisms of $\CC$, a \emph{right Bousfield localization of $\CC$ with respect to $H$} is a right $\CC$-semimodel $\XX$-category $\fromto{\CC}{R_H\CC}$ that is initial among right $\CC$-semimodel $\XX$-categories $F:\fromto{\CC}{\DD}$ with the property that for any $f$ representing a class in $H$, $\RR F(f)$ is an isomorphism of $\Ho\NN$.
\item A morphism $\fromto{A}{B}$ is a \emph{$K$-colocal equivalence} if for any representative $X$ of an element of $K$, the morphism
\begin{equation*}
\fromto{\RMor_{\CC}(X,A)}{\RMor_{\CC}(X,B)}
\end{equation*}
is an isomorphism of $\Ho s\mathrm{Set}_{\XX}$.
\item An object $Z$ of $\MM$ is \emph{$K$-colocal} if for any $K$-colocal equivalence $\fromto{A}{B}$, the morphism
\begin{equation*}
\fromto{\RMor_{\MM}(Z,B)}{\RMor_{\MM}(Z,A)}
\end{equation*}
is an isomorphism of $\Ho s\mathrm{Set}_{\XX}$.
\item A \emph{right Bousfield localization of $\CC$ with respect to $K$} is nothing more than a right Bousfield localization of $\CC$ with respect to the set of $K$-colocal equivalences.
\end{enumerate}
\end{dfn}

\begin{prp} A right Bousfield localization $R_K\CC$ is essentially unique if it exists.
\begin{proof} Initial objects are essentially unique.
\end{proof}
\end{prp}

\begin{nul} Suppose now that the set $K$ is $\XX$-small.
\end{nul}

\begin{ntn} Suppose $I$ and $J$ are generating $\XX$-small sets of cofibrations and trivial cofibrations, respectively, each with cofibrant domains. For every element $A\in K$, choose a cosimplicial resolution (i.e., a cofibrant replacement in the Reedy right semimodel category of cosimplicial objects) $\fromto{\Lambda^{\bullet}A}{A}$, and set
\begin{equation*}
I_{R_K\CC}:=J\cup\{\fromto{L_p(\Lambda^{\bullet}A)}{\Lambda^pA}\ |\ p\in\Delta, A\in K\}.
\end{equation*}
\end{ntn}

\begin{prp} \label{prp:coffibweRKC}The category $\CC$ is a structured homotopical category $R_K\CC$ with the following definitions.
\begin{enumerate}[(\ref{prp:coffibweRKC}.1)]
\item A cofibration of $R_K\CC$ is defined to be a cofibration $\fromto{X}{Y}$ of $\CC$ such that there exists a weak equivalence $\fromto{Y}{Z}$ of $\CC$ such that the composition $\fromto{X}{Z}$ is an element of $\cell I_{R_K\CC}$.
\item A fibration of $R_K\CC$ is nothing more than a fibration of $\CC$.
\item A weak equivalence of $R_K\CC$ is a $K$-colocal equivalence.
\end{enumerate}
Moreover the identity functor $\fromto{R_K\CC}{\CC}$ preserves cofibrations, where as the identity functor $\fromto{\CC}{R_K\CC}$ preserves weak equivalences.
\begin{proof} Immediate.
\end{proof}
\end{prp}

\begin{nul} Note the rough similarity between the definition of the cofibrations given here and the a priori stonger description given by Hirschhorn \cite[Proposition 5.3.6]{MR2003j:18018}, namely, than a morphism is a cofibration just in case it is a retract of an element $\fromto{X}{Y}$ of $\cell I$ such that there exists a weak equivalence $\fromto{Y}{Z}$ such that the composite $\fromto{X}{Z}$ is an element of $\cell I_{R_K\CC}$. The two are in fact equivalent when $\CC$ is a model category; this is the content of the following lemma.

If $\CC$ is not a model category, there seems to be a genuine difference between the two conditions, but, unfortunately, the distinction seems to be fairly subtle, and I do not have an example that exhibits it. In any case it is certainly the weaker of these that is needed here.
\end{nul}

\begin{lem}\label{lem:mycofRKCHirschcofRKC} Suppose $\CC$ a model category; then a morphism of $\CC$ is a cofibration of $R_K\CC$ if and only if it is a retract of an element $\fromto{X}{Y}$ of $\cell I$ such that there exists a weak equivalence $\fromto{Y}{Z}$ such that the composite $\fromto{X}{Z}$ is an element of $\cell I_{R_K\CC}$.
\begin{proof} That such a retract is a cofibration of $R_K\CC$ is obvious. This is of course true regardless of whether $\CC$ is a model category.

In the other direction, suppose $f:\fromto{X}{Y}$ a cofibration for which there is a weak equivalence $e:\fromto{Y}{Z}$ of $\CC$ such that the composite $g:\fromto{X}{Z}$ is a retract of an element of $\cell I_{R_K\CC}$. The claim is that $f$ can be written as a retract of an element $\fromto{X}{Y'}$ of $\cell I$ for which there exists a weak equivalence $\fromto{Y'}{Z}$ such that the composite $e'\circ f'=g$. Indeed, simply factor $f$, by the small object argument, as an element $f':\fromto{X}{Y'}$ of $\cell I$ followed by an element $p:\fromto{Y'}{Y}$ of $\inj I$, and set $e'=e\circ p$. Then since $\CC$ is a model category, $p$ is a trivial fibration.\footnote{If $\CC$ were not a model category, this would follow only if $Y$ were $\EE$-fibrant.} The retract argument thus implies the claim.
\end{proof}
\end{lem}

\begin{lem}\label{lem:IRKinjtrivfib} A morphism $\fromto{Y'}{Y}$ whose codomain $Y$ is fibrant is a trivial fibration in $R_K\CC$ if and only if it is an element of $\inj I_{R_K\CC}$.
\begin{proof} Of course $\fromto{Y'}{Y}$ is an element of $\inj J$ if and only if it is a fibration of $\CC$ --- and, equivalently, of $R_K\CC$. In this circumstance, for any element $A\in K$, the morphism $\fromto{\RMor_{\CC}(A,Y')}{\RMor_{\CC}(A,Y)}$ is a fibration of $s\mathrm{Set}_{\XX}$, and by \cite[Proposition 16.4.5]{MR2003j:18018}, will be an equivalence if and only if the morphism $\fromto{Y'}{Y}$ is an element of $\inj I_{R_K\CC}$.
\end{proof}
\end{lem}

\begin{lem}\label{lem:factRKCcoftrivfib} There is a functorial factorization of every morphism $\fromto{X}{Y}$ of $\CC$ with fibrant codomain $Y$ into a cofibration $\fromto{X}{Y'}$ of $R_K\CC$ followed by a trivial fibration $\fromto{Y'}{Y}$ of $R_K\CC$.
\begin{proof} The usual construction via the small object argument provides a factorization of every morphism $\fromto{X}{Y}$ (irrespective of the fibrancy of $Y$) into an element $\fromto{X}{Y'}$ of $\cell I_{R_K\CC}$ followed by an element $\fromto{Y'}{Y}$ of $\inj I_{R_K\CC}$. Now $\fromto{X}{Y'}$ is clearly a cofibration, and, by the previous lemma, $\fromto{Y'}{Y}$ is a trivial fibration of $R_K\CC$.
\end{proof}
\end{lem}

\begin{lem}\label{lem:liftcofsRKC} Cofibrations of $R_K\CC$ satisfy the left lifting property with respect to every trivial fibration of $R_K\CC$ with fibrant codomain.
\begin{proof} That this is true of any element of $I_{R_K\CC}$ is already contained in \ref{lem:IRKinjtrivfib}. It thus follows for any element of $\cell I_{R_K\CC}$.

Now suppose $\fromto{X}{Y}$ a cofibration of $\CC$ and $e:\fromto{Y}{Y'}$ a weak equivalence of $\CC$ such that the composite $\fromto{X}{Y'}$ is an element of $\cell I_{R_K\CC}$. Suppose $W$ a fibrant object, $\fromto{Z}{W}$ a trivial fibration of $R_K\CC$, and the following a commutative diagram:
\begin{equation*}
\xymatrix@C=18pt@R=18pt{
X\ar[d]\ar[r]&Z\ar[d]\\
Y\ar[d]_e\ar[r]&W\\
Y'.&
}
\end{equation*}
To prove the lemma, it will now suffice to show that there exists a lift $\fromto{Y}{Z}$. For this, note that in the right semimodel category $(X/\CC)$, $W$ is fibrant and $\fromto{Y}{Y'}$ is a weak equivalence; hence there is a homotopy lift $\fromto{Y'}{W}$ in $(X/\CC)$. Since this homotopy lift is chosen in the slice category, it follows that there is a commutative diagram
\begin{equation*}
\xymatrix@C=18pt@R=18pt{
X\ar[d]\ar[r]&Z\ar[d]\ar[d]\\
Y'\ar[r]&W.
}
\end{equation*}
Hence there is a lift $\ell:\fromto{Y'}{Z}$, and the composite $\ell\circ e:\fromto{Y}{Z}$ is thus a homotopy lift of the diagram
\begin{equation*}
\xymatrix@C=18pt@R=18pt{
&Z\ar[d]\ar[d]\\
Y\ar[r]&W
}
\end{equation*}
in $(X/\CC)$. Since $Y$ is cofibrant, the homotopy lifting property of the fibration $\fromto{Z}{W}$ implies that a strict lift of the diagram, homotopic to $\ell\circ e$, exists in $(X/\CC)$.
\end{proof}
\end{lem}

\begin{lem}\label{lem:trivcofsRKCtrivcofsC} The trivial cofibrations of $R_K\CC$ are exactly those of $\CC$.
\begin{proof} If $f:\fromto{K}{L}$ is a trivial cofibration of $\CC$, then it is a weak equivalence of $\CC$, a retract of an element of $\cell I$, and a retract of an element of $\cell J$. It follows that $f$ is a weak equivalence of $R_K\CC$ and a retract of an element of $\cell I_{R_K\CC}$ as well; thus $f$ is a trivial cofibration of $R_K\CC$.

Conversely, suppose $f:\fromto{K}{L}$ a trivial cofibration of $R_K\CC$, and let $i:\fromto{L}{L'}$ be a trivial cofibration in $\CC$ with $L'$ fibrant. Factor the composite morphism $f':=i\circ f$ as a trivial cofibration $j:\fromto{K}{K'}$ of $\CC$ followed by a fibration $p:\fromto{K'}{L'}$ of $\CC$:
\begin{equation*}
\xymatrix@C=24pt@R=24pt{
K\ar[d]_j\ar[r]^f\ar[dr]^{f'}&L\ar[d]^i\\
K'\ar[r]_p&L'.
}
\end{equation*}
Since $f'$ is a weak equivalence of $R_K\CC$, so is $p$. Now $f'$ is a trivial cofibration of $R_K\CC$ and by the previous result has the left lifting property with respect to $p$. The retract argument now implies that $f'$ is a retract of $p$ and is thus a trivial cofibration of $\CC$. Since $f'=i\circ f$, and both $f'$ and $i$ are weak equivalences of $\CC$, it follows that $f$ is a weak equivalence of $\CC$ as well. Since $f$ was a fortiori a cofibration of $\CC$, the converse is verified.
\end{proof}
\end{lem}

\begin{prp} The structured homotopical category $R_K\CC$ is an $\XX$-tractable right $\CC$-semimodel category. If in addition $\CC$ is a right proper model $\XX$-category, then $R_K\CC$ is a right proper model category (not necessarily $\XX$-tractable).
\begin{proof} The factorization axioms are \ref{lem:factRKCcoftrivfib} and the corresponding factorization in $\CC$, coupled with \ref{lem:trivcofsRKCtrivcofsC}. Likewise, the lifting properties are \ref{lem:liftcofsRKC} and the corresponding property in $\CC$, coupled with \ref{lem:trivcofsRKCtrivcofsC}.

The structured homotopical structure is by \ref{lem:mycofRKCHirschcofRKC} exactly the one provided by Hirschhorn in case $\CC$ is right proper.
\end{proof}
\end{prp}

\begin{prp} The left Quillen identity functor $U:\fromto{R_K\CC}{\CC}$ induces a coreflexive fully faithful functor of $\XX$-categories --- and thus also of $(\Ho s\mathrm{Set}_{\XX})$-categories --- $\LL U:\fromto{\Ho R_K\CC}{\Ho\CC}$. The derived right adjoint $\RR F:\fromto{\CC}{\Ho R_K\CC}$ is essentially surjective.
\begin{proof} Write $F:\fromto{\CC}{R_K\CC}$ for the left adjoint of $U$. It suffices to show that the unit $\fromto{X}{(\RR F)(\LL U)X}$ of the derived adjunction is an isomorphism of $\Ho R_K\CC$. But this is clear, as the fibrant replacement in $\CC$ of any object of $R_K\CC$ is a fibrant replacement in $R_K\CC$.
\end{proof}
\end{prp}

\begin{cor} A morphism $\fromto{A}{B}$ of $\CC$ is a $K$-colocal weak equivalence if and only if the induced morphism $\fromto{(\RR F)A}{(\RR F)B}$ is an isomorphism of $\Ho R_K\CC$.
\end{cor}

\begin{prp} A cofibrant-fibrant object $X$ of $\CC$ is cofibrant as an object of $R_K\CC$ if and only if it is $K$-colocal.
\begin{proof} Since a $K$-colocal weak equivalence $\fromto{A}{B}$ is weak equivalence of $R_K\CC$, it follows that if $X$ is cofibrant in $R_K\CC$, the morphism
\begin{equation*}
\fromto{\RMor_{\CC}(X,A)\simeq\RMor_{R_K\CC}(X,A)}{\RMor_{\CC}(X,B)\simeq\RMor_{R_K\CC}(X,B)}
\end{equation*}
is an isomorphism of $\Ho s\mathrm{Set}_{\XX}$.

On the other hand, if $X$ is $K$-colocal, it suffices, since by assumption $X$ is fibrant in $R_K\CC$, to show that the morphism $\fromto{\varnothing}{X}$ has the left lifting property with respect to all trivial fibrations $\fromto{A}{B}$ of $R_K\CC$ with fibrant codomain $B$. Since $X$ is $K$-colocal, the map
\begin{equation*}
\fromto{\Mor_{\Ho\CC}(X,A)}{\Mor_{\Ho\CC}(X,B)}
\end{equation*}
is a bijection, whence the desired lifting property.
\end{proof}
\end{prp}

\begin{cor} An object $X$ of $\CC$ is $K$-colocal if and only if the counit of the derived adjunction $\fromto{(\LL U)(\RR F)X}{X}$ is an isomorphism of $\Ho\CC$.
\end{cor}

\begin{cor} A weak equivalence of $R_K\CC$ between $K$-colocal objects is a weak equivalence of $\CC$.
\end{cor}

\begin{thm}\label{thm:rtBousfieldexists} The right $\CC$-semimodel category $R_K\CC$ is a right Bousfield localization of $\CC$ with respect to $K$.
\begin{proof} The proof that the right $\CC$-semimodel structure described here has the universal property required of a right Bousfield localization is well known, and is \cite[Proposition 3.3.18]{MR2003j:18018}.
\end{proof}
\end{thm}

\begin{cor} A cofibrant-fibrant object of $\CC$ is $K$-colocal if and only if it can be written as a homotopy colimit (in $\CC$) of a diagram of representatives of elements of $K$.
\begin{proof} The proof is in fact identical to the one given by Hirschhorn \cite[\S 5.5]{MR2003j:18018}. For any object $A$ of $\CC$, the functor
\begin{equation*}
\RMor_{\CC}(-,A):\fromto{\Ho\CC^{\op}}{\Ho(s\mathrm{Set}_{\XX})}
\end{equation*}
turns homotopy colimits into homotopy limits; hence any homotopy colimit of $K$-colocal objects is again $K$-colocal. On the other hand, one verifies easily that any $I_{R_K\CC}$-cell complex can be written as a homotopy colimit of objects of $K$. 
\end{proof}
\end{cor}

\subsection*{Application: The homotopy limit of left Quillen presheaves}{Elsewhere \cite[Theorem 2.44]{arXiv:0708.2067v1} I have constructed a model category the plays the role of the homotopy limit of a right Quillen presheaf, which I constructed by taking a left Bousfield localization of a projective model structure. Here I define a right semimodel category structure on the category of left sections of a left Quillen presheaf that plays the role of the homotopy limit; it is a right Bousfield localization of an injective model structure.}

\begin{nul} Suppose $\XX$ a universe. Suppose $K$ and $\XX$-small category, and suppose $\FF$ an $\XX$-tractable left Quillen presheaf on $K$.
\end{nul}

\begin{dfn} A left section $(X,\phi)$ of $\FF$ is said to be \emph{homotopy cartesian} if for any morphism $f:\fromto{\ell}{k}$ of $K$, the morphism
\begin{equation*}
\phi_f^h:\fromto{\LL f^{\star}X_k}{X_{\ell}}
\end{equation*}
is an isomorphism of $\Ho\FF_{\ell}$.
\end{dfn}

\begin{thm}\label{prp:holimmodcat} There exists an $\XX$-tractable right semimodel structure on the category $\Sect^L\FF$ --- the homotopy limit structure $\Sect_{\holim}^L\FF$ --- satisfying the following conditions.
\begin{enumerate}[(\ref{prp:holimmodcat}.1)]
\item\addtocounter{equation}{1} The fibrations are exactly the injective fibrations.
\item\addtocounter{equation}{1} The cofibrant objects are the injective cofibrant left sections that are homotopy cartesian.
\item\addtocounter{equation}{1} The weak equivalences between cofibrant objects are precisely the objectwise weak equivalences.
\end{enumerate}
\begin{proof} It suffices to find an $\XX$-small set $G$ of cofibrant homotopy cartesian left sections of $\FF$ such that any homotopy cartesian left section of $\FF$ can be written as a homotopy colimit of elements of $G$, for then $\Sect_{\holim}^L\FF:=R_G\Sect_{\inj}^L\FF$ satisfies the conditions of the theorem.

If $K$ is a discrete category, then the injective model category of left sections of $\FF$ is isomorphic to the product model category $\Prod_{k\in K}\FF_k$, and every left section is homotopy cartesian.

For $K$ arbitrary, the model category $\Sect_{\inj}^L\FF$ is equivalent\footnote{not merely Quillen equivalent --- equivalent as categories in a manner fully compatible with the model structure} to the injective model category of left sections of the left Quillen presheaf
\begin{equation}\label{eqn:redltQpsh}
\xymatrix@C=18pt{\Prod_{i\in\Obj K}\FF_i\ar@<0.5ex>[r]\ar@<-0.5ex>[r]&\Prod_{[j\to i]\in\nu_1K}\FF_j\ar@<1ex>[r]\ar[r]\ar@<-1ex>[r]&\Prod_{[k\to j\to i]\in\nu_2K}\FF_k},
\end{equation}
and an object of $\Sect^L\FF$ is homotopy cartesian if and only if the corresponding left section of the left Quillen presheaf \eqref{eqn:redltQpsh} is homotopy cartesian. It therefore suffices to construct such a $G$ when $K$ is of the form
\begin{equation*}
\xymatrix@C=18pt{k_2\ar@<1ex>[r]\ar[r]\ar@<-1ex>[r]&k_1\ar@<0.5ex>[r]\ar@<-0.5ex>[r]&k_0}.
\end{equation*}

For $K$ with this structures, suppose $G_0$ a set of cofibrant homotopy generators of $\FF_{k_0}$. Let $U_1$ denote the intersection of the essential images of $G_0$ under the two functors $\fromto{\Ho\FF_{k_0}}{\Ho\FF_{k_1}}$, and let $U_2$ denote the intersection of the essential images of $U_1$ under the three functors $\fromto{\Ho\FF_{k_1}}{\Ho\FF_{k_2}}$. One verifies easily that any homotopy cartesian section of $\FF$ can be written as a homotopy colimit of homotopy cartesian sections $(X_0,X_1,X_2,(\phi_f)_{f\in K})$ such that $X_0\in G_0$, $X_1\in U_1$, and $X_2\in U_2$. Let $G_1$ and $G_2$ denote $\XX$-small sets of cofibrant-fibrant representatives of the isomorphism classes of $U_1$ and $U_2$, respectively. Then the desired set $G$ of cofibrant homotopy cartesian sections is
\begin{equation*}
G:=\{(X_0,X_1,X_2,(\phi_f)_{f\in K})\ |\ X_0\in G_0,X_1\in G_1,X_2\in G_2\}.\qedhere
\end{equation*}
\end{proof}
\end{thm}

\bibliographystyle{amsplain}
\bibliography{../../math}

\end{document}